\documentclass[12pt,leqno]{amsart}
\usepackage{amsmath,amsthm,amssymb}
\usepackage{color}

\def\init{\setcounter{equation}{0}}
\newtheorem{theorem}{Theorem}[section]

\def\init{\setcounter{equation}{0}}

\newcommand{\R}{\mathbb{R}}

\newcommand{\Z}{\mathbb{Z}}

\newcommand{\T}{\mathbb{T}}
\newcommand{\red}[1]{{\color{red} #1}}

\title{
Gauge equivalence and the inverse spectral problem for the magnetic
Schr\"odinger operator on the torus
}

\author{
G. Eskin and J. Ralston,\\  Department of Mathematics, UCLA,\\
Los Angeles, CA 90095-1555, USA
}

\begin{document}
\maketitle

\begin{center}
In memory of Mark Iosifovich Vishik
\end{center}

\begin{abstract}
We study  the     inverse  spectral problem
for the Schr\"odinger operator
$H$ on the two-dimensional torus with even magnetic field $B(x)$ and
even electric potential $V(x)$.  V.Guillemin [11] proved that the spectrum
of $H$ determines $B(x)$  and $V(x)$.  A simple proof of Guillemin's results
 was given by the authors in [3]. In the present paper we consider gauge equivalent classes of magnetic potentials and give conditions which imply that the gauge
equivalence class and the spectrum of $H$  determine the magnetic field and
the electric potential. We also show that generically  the spectrum and the magnetic field determine the \lq\lq extended" gauge equivalence class of the magnetic potential.  The proof is a modification of the proof in [3]
with some corrections and clarifications.

\end{abstract}

\section{Introduction}
\init

Let
$L=\{m_1e_1+m_2e_2: m=(m_1,m_2)\in \Z^2\}$
be  a lattice in $\R^2$.  Here  $\{e_1,e_2\}$  is a basis in $\R^2$.
We assume that the lattice $L$  has the following property:
\begin{equation}                                      \label{eq:1.1}
\mbox{For}\ \ d,d'\in L,\ \mbox{if}\ |d|=|d'|,\ \ \mbox{then}\ \ d'=\pm d.
\end{equation}
Let $L^*=\{\delta\in\R^2: \delta\cdot d\in \Z \ \mbox{for all}\ \ d\in L\}$ be
the dual lattice.  We consider
a Schr\"odinger operator  of the form
\begin{equation}                                    \label{eq:1.2}
H=\Big(-i\frac{\partial}{\partial x_1}-A_1(x)\Big)^2
+\Big(-i\frac{\partial}{\partial x_1}-A_2(x)\Big)^2 +V(x),\ \ x\in \R^2,
\end{equation}
where
$A(x)=(A_1(x),A_2(x))$  is the magnetic potential and $V(x)$  is the electric
potential.

Let $B(x)$  be the magnetic field,
\begin{equation}                                   \label{eq:1.3}
B(x)= \mbox{curl}\, A(x)=\frac{\partial A_2}{\partial x_1}
-\frac{\partial A_1}{\partial x_2}.
\end{equation}
We assume that $B(x)$ and $V(x)$ are periodic,  i.e.
$$
B(x+d)=B(x),\ \ V(x+d)=V(x) ,\ \ \forall d\in L,
$$
i.e.  $B(x)$  and $V(x)$  are smooth functions  on $\T^2=\R^2/L$.
We also assume that $B(x)$  and $V(x)$  are even,  i.e.
$B(-x)=B(x)$  and $\ V(-x)=V(x)$.

Denote by $G(\T^2)$
the gauge group of complex-valued functions $g(x)\in C^\infty(\T^2)$
such that
$|g(x)|=1$.
Any $g(x)\in G(\T^2)$  has the following form
\begin{equation}                                \label{eq:1.4}
g(x)=\exp(2\pi i\delta\cdot x + i\varphi(x)),
\end{equation}
where $\delta\in L^*$  and $\varphi(x)$  is periodic,
$\varphi(x+d)=\varphi(x),\ \ \forall d \in L$.
The operator of multiplication by $g(x)$ transforms
the equation $Hu=\lambda u$  to the equation
$$
H' u'(x)=\lambda u'(x),
$$
where
$H'$ has the form (\ref{eq:1.2}) with $A(x)$ replaced by $A'(x)$,
\begin{equation}                                  \label{eq:1.5}
A'(x)=A(x)-ig^{-1}(x)\nabla g(x)=A(x)+2\pi\delta +\nabla\varphi(x),\ \delta\in L^*,
\end{equation}
  and $u'(x)=g^{-1}(x)u(x)$.
The magnetic potentials $A'(x)$  and $A(x)$
related by (\ref{eq:1.5}) are called gauge equivalent.
Since $H$  and $H'$  are unitarily equivalent, they have the same spectrum.

Let
\begin{equation}                                      \label{eq:1.6}
B(x)=\sum_{\beta\in L^*}b_\beta e^{2\pi i\beta\cdot x}
\end{equation}
be
the Fourier series expansion of $B(x)$.
We assume that the coefficient
\begin{equation}                               \label{eq:1.7}
b_0=|D|^{-1}\int\limits_{D}B(x)dx,
\end{equation}
 is not zero.
 Here
$D$  is a fundamental domain for the lattice $L$ given by
\begin{equation}                             \label{eq:1.8}
D=\{t_1e_1+t_2e_2,\ |t_j|\leq \frac{1}{2}, j=1,2\}
\end{equation}
  and
$|D|$  is the area of $D$.  Note that if $x\in D$,  then $-x\in D$.

Given $B(x)$,  we let
\begin{equation}                                    \label{eq:1.9}
A(x)=A^{0}(x)+a_0+\sum_{\beta\in (L^*\setminus 0)}a_\beta e^{2\pi i\beta\cdot x},
\end{equation}
where $a_0=(a_{01},a_{02})$ is a constant,
\begin{equation}                             \label{eq:1.10}
A^0(x)=\frac{b_0}{2}(-x_2,x_1)
\end{equation}
and
\begin{equation}                              \label{eq:1.11}
a_\beta=b_\beta(2\pi i)^{-1}(\beta_1^2+\beta_2^2)^{-1}(-\beta_2,\beta_1).
\end{equation}
Note that $\mbox{curl}\,A=B(x)$.

Let $A'(x)$  be any magnetic potential that is gauge
equivalent to $A(x)$.
Since $\mbox{curl}\,\nabla\varphi =0$\red{,}
$$
B'(x)=B(x),
$$
where $B'(x)=\mbox{curl}\, A'(x)$.  Therefore $A'(x)$  can be also
represented in the form (\ref{eq:1.9}) with $a_0^{^\prime}$ not necessarily equal to $a_0$.
 For the gauge equivalence of $A^\prime(x)$ and $A(x)$, in addition to the equality of the magnetic fields, one needs (cf. (\ref{eq:1.5}))

\begin{equation}                                      \label{eq:1.12}
a_0'-a_0=2\pi\delta,
\end{equation}
 for some $\delta\in L^*.$
The main question in this paper is: For the class of magnetic potentials considered here, to what extent do the spectra of magnetic Schr\"odinger
operators and the gauge equivalence classes of magnetic potentials determine the magnetic and electric fields?

In \S 2  we will describe  the domain on which $H$  is a self-adjoint
operator  with compact resolvent and hence  has a discrete spectrum.
Here we will describe the gauge equivalence classes of the magnetic potential assuming
that the magnetic field $B(x)$  is fixed.

Let $\gamma_j,j=1,2,$  be the basis of the homology group of the torus,
given by
$\gamma_j=\{te_j,0\leq t\leq 1\}$.  Let
$$
\alpha_j=\int\limits_{\gamma_j} a_0\cdot dx=a_0\cdot e_j,
$$
where
$a_0$  is  the constant vector in (\ref{eq:1.9}).  For any  $d=m_1e_1+m_2e_2\in L$
we have $a_0\cdot d= m_1\alpha_1+m_2\alpha_2$, i.e. knowing $\{\alpha_1,\alpha_2\}$
determines  $a_0\cdot d$  for any  $d\in L$.

Let $A'(x)$  be a magnetic potential of the form (\ref{eq:1.9})
with $a_0$  replaced by $a_0'$.  Define
$$
\alpha_j'=
\int\limits_{\gamma_j}a_0'\cdot dx=a_0'\cdot e_j.
$$
Let $\{e_1^*,e_2^*\}$  be the basis in $L^*$  dual to $\{e_1,e_2\}$,
i.e. $e_j\cdot e_k^*=\delta_{jk}$.

 The potentials $A(x)$  and $A'(x)$  are gauge equivalent if and only if  $\mbox{curl}\,A =\mbox{curl}\, A'$ and (cf.  (\ref{eq:1.12}))
\begin{equation}                                   \label{eq:1.13}
\alpha_j'-\alpha_j=(a^\prime_0-a_0)\cdot e_j=2\pi\delta\cdot e_j,\ j=1,2,
\end{equation}
for some
$\delta\in L^*.$
Short equivalent
 forms of (\ref{eq:1.13}) are
$$
e^{i\alpha_j}=e^{i\alpha_j^\prime},\ j=1,2,
$$
and
\begin{equation}                                   \label{eq:1.14}
e^{ia_0\cdot d}=e^{ia_0'\cdot d}\ \ \mbox{for any}\ \ d=n_1e_1+n_2e_2\in L.
\end{equation}

Changing $x$  to $-x$  we get  the operator $H'$  which is just $H$ with $a_0$ changed to $a_0^\prime=-a_0$.
Note that  $H$  and $H'$ have the same spectrum
 but their magnetic potentials are not gauge equivalent when $a_0\neq 0$. Since we are looking for consequences of isospectrality, we introduce a weaker notion of gauge equivalence, namely
\begin{equation}                           \label{eq:1.15}
\cos a_0\cdot d=\cos a_0'\cdot d,\ \forall d\in L
\end{equation}
The condition (\ref{eq:1.15}) is equivalent to $\cos\alpha_j=\cos\alpha_j^\prime,\ j=1,2$.
 Since
$$\cos\alpha_j-\cos\alpha_j^\prime=2\sin(\frac{\alpha_j-\alpha_j^\prime}{2})\sin(\frac{\alpha_j+\alpha_j^\prime}{2}),$$
$\cos\alpha_j=\cos\alpha_j^\prime$ implies that either $\alpha_j-\alpha_j^\prime$ or $\alpha_j-\alpha_j^\prime$ is an integer multiple of $2\pi$. Thus there are two choices for each $j$: $e^{ia_0\cdot e_j}=e^{ia_0^\prime\cdot e_j}$ or $e^{ia_0\cdot e_j}=e^{-ia_0\cdot e_j}$. We will say that $a_0'$  and $a_0$  belong to
the same ``extended gauge equivalence class"
if (\ref{eq:1.15})  holds. Thus for every extended gauge equivalence class of magnetic potentials, there are four choices of $a_0$, including $a_0^\prime =a_0$ and $a_0^\prime =-a_0$, giving distinct gauge equivalence classes when $a_0\neq 0$.
\vskip.2in
 Our first result gives conditions for the spectrum of $H$ and gauge equivalence class of $A$ to determine the fields:
\begin{theorem}                                   \label{theo:1.1}
Let $B(x),V(x)$ be periodic and even smooth functions,
and assume $L$  satisfies the condition (\ref{eq:1.1}).
Suppose
\begin{equation}                                  \label{eq:1.16}
\int\limits_D B(x)dx=2\pi.
\end{equation}
Consider the spectrum of the Schr\"odinger operator $H$
with $A(x)$ having the form (\ref{eq:1.9}).
Suppose that
\begin{equation}                           \label{eq:1.17}
|B(x)-b_0|<|b_0|,
\end{equation}
where $b_0=2\pi/|D|$.
Then the spectrum of $H$  determines uniquely $B(x)$  and  $V(x)$
assuming that $\cos  a_0\cdot d,\ d\in L,$
is given and
\begin{equation}                                 \label{eq:1.18}
\cos a_0\cdot d\neq 0\ \ \mbox{for all}\ \ d\in L.
\end{equation}
\end{theorem}
\vskip.2in
\begin{theorem}                                \label{theo:1.2}
 Assume that  the conditions (\ref{eq:1.1}), (\ref{eq:1.16}), (\ref{eq:1.17})
hold, and that the spectrum  of $H$ and the magnetic field $B(x)$ are given.
If $B(x)$ satisfies a generic condition (stated in the proof),
then $\cos a_0\cdot d$  is determined for all $d\in L$,  i.e.  the exended gauge equivalence
class of $A(x)$  is determined by the spectrum and the magnetic field.
\end{theorem}
It follows from  Theorem \ref{theo:1.2}   that
if $B(x),V(x)$  are fixed  and  the extended gauge equivalence classes of
$A(x)$  and $A'(x)$ are different,
i.e. if $\cos a_0\cdot d\neq \cos a_0'\cdot d$ for some $d\in L$,
  then the corresponding
operators $H$  and $H'$  have different spectra.
This confirms  the Aharonov-Bohm effect  stating  that
 different gauge equivalence classes have a different quantum mechanical
effects,  for example,  the corresponding
 Schr\"odinger  operators have different spectra(cf. [5]).

The case when $B(x)$  and $V(x)$ are even  and $a_0=0$ (cf. (\ref{eq:1.9}))
was proven in an important paper of Guillemin [11]. In [3] we reproved the result
of [6] by a different and simpler method.

The method of the present paper is a modification of the method of [3] with some
clarifications  and  corrections.

We mention briefly some related results on the inverse spectral problems
in two and higher dimensions:
The magnetic Schr\"odinger operator on $\T^2$
with $\int_D B(x)dx=0$ and $A(x)$ periodic was studied in [1].   In [10] Gordon et al. generalized [11] to the case of
$n$-dimensional tori.
Guillemin and Kazhdan [13]  studied  the inverse  spectral  problem for
negatively curved manifolds.  Guillemin [12] studied the inverse
spectral problem on $S^2$.  Zeldich [18] solved the inverse spectral  problem for
analytic bi-axisymmetric plane  domains.  In [6]  the inverse spectral problems
on the torus for the Schr\"odinger operator $-\Delta+q(x)$  were studied.  See
also [4],[8], [16].

Gordon [7] and Gordon-Schuth [9]  gave many interesting examples of
isospectral manifolds { which were not isometric.

\section{The singularities of the wave trace}
\init

We introduce the ``magnetic translation operators" (cf. [17])
\begin{equation}                               \label{eq:2.1}
T_ju(x)=e^{-iA^0(e_j)\cdot x}u(x+e_j),\ \ j=1,2,
\end{equation}
where
$A^0(x)$  is from (\ref{eq:1.10}).
 These operators are required to commute with each other and with $H$. This implies that
\begin{equation}                                      \label{eq:2.2}
A^0(e_1)\cdot e_2=-A^0(e_2)\cdot e_1=\pi l,
\end{equation}where $l$  is an integer.  Using
(\ref{eq:1.7}),  (\ref{eq:1.10})  we get that (\ref{eq:2.2})
is equivalent  to
\begin{equation}                                      \label{eq:2.3}
\int\limits_D B(x)dx_1dx_2=2\pi l.
\end{equation}
Later we shall assume that $l=1$.
Having that $T_1,T_2$ and $H$ commute we denote by $D_0$  the subspace
of the Sobolev space $H^2(\R^2)$ consisting of $u(x)\in H^2(\R^2)$  such
that $T_ju=u,j=1,2$.
Then  the operator $H$ is self-adjoint in $L_2(D)$  on the restriction of $D_0$
to the
fundamental domain $D$.
We shall denote this operator  by $H_D$.

Let $\lambda_1\leq\lambda_2\leq \lambda_3\leq ...$  be the spectrum  of
$H_D$.   and let $E_D(x,y,t)$ be the  fundamental solution for the wave equation on $\Bbb R^2/L$. Then the
 wave trace formula   gives the equality as distributions in $t$
\begin{equation}                              \label{eq:2.5}
\sum_{j=1}^\infty \cos t\sqrt{\lambda_j}=\int\limits_D E_D(x,x,t)dx.
\end{equation}
The distribution $E_D(x,y,t)$ is defined as follows:
Let $E(x,y,t)$  be the fundamental solution for  the wave equation on $\Bbb R^2$:
\begin{align}                              \label{eq:2.6}
&\frac{\partial^2 E(x,y,t)}{\partial t^2}+HE(x,y,t)=0,\ \ x\in \R^2,y\in \R^2,
\\
\nonumber
& E(x,y,0)=0,\ \ \frac{\partial E(x,y,0)}{\partial t} =\delta(x-y),
\ \ x\in \R^2,y\in \R^2.
\end{align}
Then
\begin{equation}                                 \label{eq:2.7}
E_D(x,y,t)=\sum_{(m,n)\in \Z^2}T_1^mT_2^nE(x,y,t).
\end{equation}
Note that
\begin{equation}                                \label{eq:2.8}
T_1^mT_2^nE(x,y,t)=e^{-iA^0(d)\cdot x}E(x+d,y,t),
\end{equation}
where
$d=me_1+ne_2$.  We used in (\ref{eq:2.8})  that $A^0(e_j)\cdot e_j=0$.
Since $E(x,y,t)$ is singular  only  when $|x-y|^2=t^2$  and
since  the condition (\ref{eq:1.1}) holds,  the singularities of the trace
(\ref{eq:2.5})   at $t=|d|,\ d=m_1e_1+m_2e_2$
come only from two terms
\begin{equation}                               \label{eq:2.9}
\int\limits_D(T_1^mT_2^nE(x,x,t)+T_1^{-m}T_2^{-n}E(x,x,t))dx.
\end{equation}
To compute the singularities in (\ref{eq:2.9})  we will use
as in [3] and [4] the Hadamard-H\"ormander
parametrix (cf. [14], [15]).  We have
\begin{equation}                                \label{eq:2.10}
E(x,y,t)=\frac{\partial}{\partial t}(E_+(x,y,t)-E_+(x,y,-t),
\end{equation}
where $E_+(x,y,t)$  is  the forward fundamental solution:
\begin{equation}                                  \label{eq:2.11}
\Big(\frac{\partial^2}{\partial t^2}+H\Big)E_+(x,y,t)=\delta(t)\delta(x-y),
\ \ E_+(x,y,t)=0\ \mbox{for}\ t<0.
\end{equation}
It follows from [14]  that
\begin{multline}                                    \label{eq:2.12}
E_+(x,y,t)=m_0(x,y)\frac{1}{2\pi}(t^2-|x-y|^2)_+^{-\frac{1}{2}}
\\
+m_1(x,y)2^{-2}\pi^{-1}(t^2-(x-y)^2)_+^{\frac{1}{2}}
+ O((t^2-(x-y)^2)^{\frac{3}{2}}),
\end{multline}
where
\begin{equation}                                   \label{eq:2.13}
m_0(x,y)=\exp\Big(i\int\limits_0^1(x-y)\cdot A(y+s(x-y))ds\Big),
\end{equation}
\begin{equation}                                   \label{eq:2.14}
m_1(x,y)=-m_0(x,y)\Big(\int\limits_0^1V(y+s(x-y))ds+b(x,y)\Big),
\end{equation}
where
$$
(2.14')\ \ \ \ \ \ \ \
b(x,y)=\Big[\Big(-i\frac{\partial}{\partial x}-A(x)\Big)^2m_0\Big]m_0^{-1}(x,y).
\ \ \ \ \ \ \ \ \ \ \ \ \ \ \ \ \ \ \ \ \ \ \
$$
Let
\begin{equation}                                         \label{eq:2.15}
I(d)=\int\limits_D\exp i\Big[-A^0(d)\cdot x+\int\limits_0^1 (A(x+sd)\cdot d) ds\Big]dx
\end{equation}
\begin{equation}                                          \label{eq:2.16}
J(d)=\int\limits_D\int\limits_0^1[(V(x+sd)+b(x+sd,x)]ds \ \exp\Big[
-iA^0(d)\cdot x+i\int\limits_0^1(A(x+sd)\cdot d)ds\Big]dx
\end{equation}
It follows from (\ref{eq:2.9})-(\ref{eq:2.14})  that  $I(d)+I(-d)$
and $J(d)+J(-d)$ are determined by the spectrum  of $H$.

Using that $A^0(d)\cdot d=0$  and that $A^0(d)\cdot x=-A^0(x)\cdot d$
we can rewrite (\ref{eq:2.15})  in the form
\begin{equation}                                 \label{eq:2.17}
I(d)=\int\limits_D \exp i\Big[2A^0(x)\cdot d+a_0\cdot d+
\int\limits_0^1(A^1(x+sd)\cdot d) ds\Big]dx,
\end{equation}
where
\begin{equation}                                  \label{eq:2.18}
A^1(x)=\sum_{\beta\in L^*\setminus 0}a_\beta e^{2\pi i\beta\cdot x},
\end{equation}
and
$a_\beta$ are defined in (\ref{eq:1.11}).

Let $\{e_1^*,e_2^*\}$  be the basis in $L^*$  dual to the basis $\{e_1,e_2\}$, i.e.
\begin{equation}                                   \label{eq:2.19}
e_j^*\cdot e_k=\delta_{jk}, \ \  1\leq j,k\leq 2.
\end{equation}
We shall construct $e_j^*,j=1,2,$  explicitly. Let $e_j^\perp=(-e_{j2},e_{j1})$,
Denote by $\Delta$  the determinant
$\begin{vmatrix}
e_{11} &e_{12}\\
e_{21} &e_{22}
\end{vmatrix}.
$
We assume that $\Delta>0$.  Note that $\Delta$ is the area of the fundamental domain
$D: \Delta=|D|$.  Now we define
\begin{equation}                                      \label{eq:2.20}
e_1^*=-\frac{1}{\Delta} e_2^\perp,\ \ \ e_2^*=\frac{1}{\Delta}e_1^\perp.
\end{equation}
We have $d=me_1+ne_2=k(m_0e_1+n_0e_2)$,  where $k\geq 1$,
is an integer and $m_o,n_0$  have no common factors.  Then
$$
A^0(d)=\frac{b_0}{2}d^\perp =\frac{b_0k}{2}(m_0e_1^\perp+n_0e_2^\perp)
=\frac{b_0k\Delta}{2}(m_0e_2^*-n_0e_1^*)=\frac{b_0k\Delta}{2}\delta,
$$
where $\delta=-n_0e_1^*+m_0e_2^*$.

Using
(\ref{eq:1.7}), (\ref{eq:2.3})  we get
\begin{equation}                            \label{eq:2.21}
A^0(d)=\pi kl\delta.
\end{equation}
If $\beta\cdot d\neq 0$,  then $\int_0^1e^{2\pi is\beta\cdot d}ds=0$.
When $\beta\cdot d=0$, we have $\beta=p\delta,\ p\in \Z\setminus O,\
\delta=-n_0e_1^*+m_0e_2^*$.  We shall compute the inner product
$d\cdot a_{p\delta}$,  where
$a_\beta$  is given by (\ref{eq:1.11}).
Note that $d=k(m_0e_1+n_0e_2),\ \delta^\perp=-n_0(e_1^*)^\perp+m_0(e_2^*)^\perp=
\frac{1}{\Delta}(n_0e_2+m_0e_1)$  (cf. (\ref{eq:2.20}).  Therefore
$$
d\cdot a_{p\delta}=\frac{kb_{p\delta}p |m_0e_1+n_0e_2|^2}
{\Delta 2\pi i p^2|\delta_\perp|^2}
=
\frac{kb_{p\delta}\Delta}{2\pi i p},
$$
since
$|\delta^\perp|^2=\frac{1}{\Delta^2}|m_0e_1+n_0e_2|^2$.  It follows
from (\ref{eq:1.7})
and (\ref{eq:2.3}) that
$\Delta=|D|=\frac{2\pi l}{b_0}$.  Hence
\begin{equation}                           \label{eq:2.22}
d\cdot a_{p\delta}=\frac{klb_{p\delta}}{ip b_0}.
\end{equation}
Therefore $I(d)$  has now the form
\begin{equation}                                  \label{eq:2.23}
I(d)=
e^{ia_0\cdot d}\int\limits_D\exp 2\pi i kl((x\cdot\delta)+A_\delta^1(\delta\cdot x)]dx,
\end{equation}
where
\begin{equation}                                 \label{eq:2.24}
A_\delta^1(\delta\cdot x)=
\sum_{p\neq 0}\frac{b_{p\delta}}{2\pi ipb_0}e^{2\pi i p(\delta\cdot x)}.
\end{equation}
 Potentials of the form (\ref{eq:2.24}) are called `` directional potentials"
in [6].  Note that
$$
\frac{d}{ds}A_\delta^1(s)=\frac{1}{b_0}B_\delta(s),
$$
where
\begin{equation}                              \label{eq:2.25}
B_\delta(s)=\sum_{p\neq 0}b_{p\delta}e^{ips}
\end{equation}
is a directional potential for $B(x)$.
>From here on we assume that $l=1$ and set $d_0=\frac{1}{k}d=m_0e_1+n_0e_2$.

Choose $\delta'\in L^*$  so that
$(\delta,\delta')$  is a basis in $L^*$
 and let $(\gamma,\gamma')$  be
the basis in $L$ dual  to $(\delta,\delta')\in L^*.$
 We let  $D'$  be the fundamental domain for $L$  with respect
 to the  basis $\{\gamma,\gamma'\}$ in
 the form
$\{s\gamma+s'\gamma',-\frac{1}{2}\leq s,s'\leq \frac{1}{2}\}$, and continue to let  $D$  be the fundamental domain  for $L$ from (\ref{eq:1.8}).
Setting $x=s\gamma+s'\gamma'$, we $x\cdot\delta=s,
x\cdot\delta'=s'$.  Since the image of $D$ is $D'$,  using this change of variables we get
\begin{equation}                                   \label{eq:2.26}
\int\limits_D\exp 2\pi ik((x\cdot\delta)+A_\delta^1(x\cdot \delta))dx=
\int\limits_{-1/2}^{1/2}
\,\int\limits_{-1/2}^{1/2}
[\exp 2\pi i k (s+A_\delta^1(s))]c_0dsds',
\end{equation}
where $c_0=\big|\frac{\partial(x_1,x_2)}{\partial(s,s')}\big|$  is
the Jacobian.
Note  that $c_0$ is the area of  a fundamental domain for $L$.
Therefore,  integrating in $s'$ we get
$$
I(d)=c_0\int\limits_{-1/2}^{1/2}
\exp \Big[2\pi i k \Big(s+\frac{a_0\cdot d_0}{2\pi}+A_\delta^1(s)\Big)\Big]ds.
$$
When  $a_0=0$  and $A_\delta^1(s)$  is an odd function of $s$  then $I(-d)=I(d)$
and the computations are simplified.  When $a_0\neq 0$   the spectral
invariant  is $I(d)+I(-d)$  and we get, changing $k$ to $-k$:
\begin{equation}                                        \label{eq:2.27}
I(d)+I(-d)=2c_0
\int\limits_{-1/2}^{1/2}
\cos \Big[2\pi  k \Big(s+\frac{a_0\cdot d_0}{2\pi}+A_\delta^1(s)\Big)\Big]ds
\end{equation}

Since $\int\limits_{-1/2}^{1/2}
e^{2\pi i(x+sd)\cdot\delta'}ds=0$ if $d\cdot\delta'\neq 0$  and it is equal $e^{2\pi i x\cdot\delta}$
 when $\delta\cdot d=0$,
the directional potential $B_\delta$ is given by
$$B_\delta(x)=\int\limits_{-1/2}^{1/2}
(B(x+sd)-b_0)ds.$$
Hence, by the assumption $|B(x)-b_0|<|b_0|$
\begin{equation}                                  \label{eq:2.28}
\max |B_\delta(x)|\leq \int\limits_{-1/2}^{1/2}
\max|B(x)-b_0|ds<|b_0|.
\end{equation}
Letting
\begin{equation}                                  \label{eq:2.29}
y=s+A_\delta^1(s),
\end{equation}
we have
$$
\frac{dy}{ds}=1+\frac{d}{ds}A_\delta^1(x)=1+\frac{B_\delta(s)}{b_0}>0.
$$
  Therefore
 the inverse function $s=s(y),\ y\in\R^1$,  is defined.  Since $y=y(s)$  is odd,
the  inverse function $s=s(y)$  is also odd.  Since $y(s+1)=y(s)+1$  we have
$s(y+1)=s(y)+1$.  Differentiating in $y$  we get $s'(y+1)=s'(y),$
i.e.  $s'(y)$  is periodic of period 1.   The function $s'(y)$  is even  since
$s(y)$  is odd.  Let  $e(y)=s(y)-y$.  Since $s(y+1)=s(y)+1$  we get   that
$e(y+1)=e(y),$  i.e.  $e(y)$  is periodic  and odd.  Note that
$$
s'(y)=1+e'(y).
$$
After
the change of variables  $s=s(y)$ in (\ref{eq:2.27}), we have
$$I(d)+I(-d)=2c_0\int\limits_{-1/2+A_\delta^1(-1/2)}^{1/2+
A_\delta^1(1/2)}
\cos (2\pi  k y+ka_0\cdot d_0)s'(y)dy.
$$
Since $\cos (2\pi ky+ka_0\cdot d_0)s'(y)$  is periodic and
$A_\delta^1(1/2)=A_\delta^1(-1/2)$  we get
\begin{equation}                         \label{eq:2.30}
I(d)+I(-d)=2c_0\int\limits_{-1/2}^{1/2}
\cos ( 2\pi  k y+ka_0\cdot d_0))s'(y)dy.
\end{equation}
Since $(\sin 2\pi ky)s'(y)$  is an odd function on $(-1/2,1/2)$,$$
2c_0\int\limits_{-1/2}^{1/2}
(\sin  ka_0\cdot d_0)(\sin  2\pi  k y) s'(y)dy=0,
$$ and this implies
\begin{equation}                         \label{eq:2.31}
I(d)+I(-d)=2c_0\int\limits_{-1/2}^{1/2}
(\cos  ka_0\cdot d_0)(\cos  2\pi  k y) s'(y)dy
\end{equation}

As an even smooth function, $s'(y)$  is a sum of  its Fourier cosine series on
$(-1/2,1/2)$.  Suppose $\cos  (ka_0\cdot d_0)$  is known and
nonzero for all $k\geq 1$.  Then we know the  Fourier cosine coefficients
for $k\geq 1$.
This uniquely determines $s'(y)$  up to a constant.  Therefore
$s'(y)=C+s_1(y)$,  where
$\int\limits_{-1/2}^{1/2} s_1(y)dy=0$.
Since
$s(y)=y+e(y)$,  where $e(y)$  is periodic and odd we get
$s'(y)=1+e'(y)$.
Knowing $s'(y)$ we can find $s(y)$  and subsequently $A_\delta^1(s)$  from
the knowledge of the spectrum of $H$  and $\cos(ka_0\cdot d_0),\ k\geq 1$.
Repeating the same arguments for any $d\in L$ we can recover $A^1(x)$.
\vskip.2in
Now we can recover $V(x)$  assuming that $\int_D V(x)dx=0$.
One can check from (2.14$'$) that
$b(x,y)$  does not depend on $a_0$:
the only term in $m_0(x,y)$  (cf.  (\ref{eq:2.13})) that contains $a_0$
has
the form $e^{ia_0\cdot(x-y)}$.  Therefore
$-i\frac{\partial}{\partial x}m_0(x,y)=(a_0+c_1(x,y))m_0(x,y)$  where $c_1(x,y)$  is independent of $a_0$.  Hence  $(-i\frac{\partial}{\partial x}-A)m_0(x,y)=(a_0+c_1(x,y)-
(a_0+\frac{b_0}{2}x^\perp+A^1(x))m_0(x,y)$  will not contain $a_0$.
  Therefore $b(x,y)$  is
known once we know $A^1(x)$.

Since $\int_0^1V(x+sd)ds=V_\delta(x)$,  we have
\begin{equation}                            \label{eq:2.32}
J(d)=J_1(d)+J_2(d),
\end{equation}
where
\begin{equation}                             \label{eq:2.33}
J_1(d)=\int\limits_D V_\delta(x\cdot \delta)\exp [2\pi i k(x\cdot\delta +\frac{a_0\cdot d_0}{2\pi}+
A_\delta^1(x\cdot \delta))]dx,
\end{equation}
and
$J_2(d)$  is the term containing $b(x,y)$,  i.e.
$J_2(d)$  is known.
Therefore $J_1(d)+J_1(-d)$  is determined by the spectrum assuming that
$\cos ka_0\cdot d$  is known.

Making the change of variables $x=s\gamma+s'\gamma'$  as in (\ref{eq:2.26})  and
integrating in $s'$  we get
\begin{equation}                                \label{eq:2.34}
J_1(d)+J_1(-d)=
2c_0\int\limits_{-1/2}^{1/2}
V_\delta(s)\cos 2\pi k (s+\frac{a_0\cdot d_0}{2\pi}+A_\delta^1(s))ds.
\end{equation}
Making the change of variables $y=s+A_\delta^1(s)$ as in (\ref{eq:2.30})
we get
\begin{equation}                                \label{eq:2.35}
J_1(d)+J_1(-d)=
2c_0\int\limits_{-1/2}^{1/2}
V_\delta(s(y))s'(y)\cos (2\pi  ky +ka_0\cdot d_0)dy,
\end{equation}
where $s=s(y)$  is the inverse to $y=s+A_\delta^1(s)$.
Note that $s'(y)$  is even periodic function of period 1,
and $V_\delta(s(y))$ is also even periodic,
since $V(x)$  is even periodic  and $s(y)$  is an odd function  satisfying $s(y+1)=s(y)+1$.

Since $V_\delta(s(y))s'(y)$
is an even function, $ \int\limits_{-1}^{1}
V_\delta(s(y))s'(y)\sin  2\pi  k ydy=0.$
Thus, as in (\ref{eq:2.31})  we have:
\begin{equation}                                \label{eq:2.36}
J_1(d)+J_1(-d)=
2c_0
\cos (ka_0\cdot d_0)
\int\limits_{-1/2}^{1/2}
V_\delta(s(y))s'(y)\cos  2\pi  k ydy, \ \ k\geq 1.
\end{equation}

Knowing
$J_1(d)+J_1(-d)$  and $\cos (ka_0\cdot d_0)$  we
know the Fourier cosine coefficients  of  the even function $V_\delta(s(y))s'(y)$  for $k\geq 1$.
Therefore we can determine $V_\delta(s(y))s'(y)$  up to a constant
\begin{equation}                                     \label{eq:2.37}
V_\delta(s(y))s'(y)=C+e_1(y),
\end{equation}
where $e_1(y)$  is known and $\int_{-1/2}^{1/2}
e_1(y)dy =0$.

Integrating (\ref{eq:2.37}) we get
$$
C=\int\limits_{-1/2}^{1/2}
V_\delta(s(y))s'(y)dy
=\int\limits_{-1/2+e(-1/2)}^{1/2+
e(1/2)}V_\delta(s)ds
=
\int\limits_{-1/2}^{1/2}
V_\delta(s)ds=0,
$$
where $s(y)=y+e(y),$  and $e(y)$  is periodic with period 1.
Thus
we can recover $V_\delta(s)$  for each $\delta\in L^*$  and therefore
we can recover $V(x)$.
This concludes the proof of Theorem \ref{theo:1.1}.

Now we shall prove  Theorem \ref{theo:1.2}.
We shall assume that the magnetic field is generic  in the following sense:
There are  two directions  $\delta_1$  and
  $\delta_2$  which form a basis for
$L^*$  such that  the directional
fields $B_{\delta_1}(s)$  and $B_{\delta_2}(s)$ are not identically zero.
In this case the functions $s_1(x)$  and $s_2(x)$ (cf. (\ref{eq:2.29}))
corresponding to $\delta_1$  and $\delta_2$
respectively,  are not identically  zero. We make the additional generic assumption that
$a_{1j}\neq 0,j=1,2,$
where $a_{kj}$ are the Fourier cosine coefficients of $s_j'(y),k\geq 0$.
Then from the main relation (\ref{eq:2.31}) we can recover
$\cos a_0\cdot d_j,j=1,2,$
where $\{d_1,d_2\}$ is the dual basis to $\{\delta_1,\delta_2\}$. Given $d\in L$, there are are integers $m$ and $n$ such that $d=md_1+nd_2$.  Hence
$$\cos a_0\cdot d=Re\{e^{ia_0\cdot(md_1+nd_2)}\}=$$ $$ Re\{(\cos(a_0\cdot d_1)\pm i\sqrt{1-(\cos(a_0\cdot d_1)^2})^m (\cos(a_0\cdot d_2)\pm i\sqrt{1-(\cos(a_0\cdot d_2)^2})^n\},$$  and, since the $\pm$'s disappear when one takes the real part, $\cos (a_0\cdot d)$ is determined by $\cos (a_0\cdot d_1)$ and $\cos (a_0\cdot d_2)$. Thus,  $\cos a_0\cdot d=\cos a_0'\cdot d$  for all
$d\in L$ as in (cf. (\ref{eq:1.15})). This proves Theorem \ref{theo:1.2}.

\bibliographystyle{amsalpha}

\end{document}